\documentclass[titlepage,11pt]{article}
\usepackage{latexsym}
\usepackage{amsfonts}
\usepackage{amsmath}
\newcommand\blackslug{\hbox{\hskip 1pt \vrule width 4pt height 8pt depth 1.5pt
        \hskip 1pt}}
\newcommand\bbox{\hfill \quad \blackslug \bigbreak}

\def\ll{,\ldots,}

\title{Pure pairs. II. Excluding all subdivisions of a graph}
\author{Maria Chudnovsky\thanks{Supported by NSF grant DMS-1550991.
This material is based upon work supported in part by the U. S. Army
Research Laboratory and the U. S. Army Research Office under grant
number
W911NF-16-1-0404.}\\
Princeton University, Princeton, NJ 08544
\\
\\
Alex Scott\thanks{Supported by a Leverhulme Trust Research
Fellowship}\\
Mathematical Institute, University of Oxford, Oxford OX2 6GG, UK
\\
\\
Paul Seymour\thanks{Supported by ONR grant N00014-14-1-0084 and NSF
grant DMS-1265563.} and Sophie Spirkl\\
Princeton University, Princeton, NJ 08544}

\date{October 28, 2017; revised \today}

\newtheorem{thm}{}[section]
\newcommand{\Proof}{\noindent{\bf Proof.}\ \ }

\begin{document}
\maketitle
\begin{abstract}
We prove for every graph $H$ there exists $\epsilon>0$ such that, for every graph $G$ with $|G|\ge 2$, 
if no induced subgraph of $G$
is a subdivision of $H$, then either some vertex of $G$ has at least $\epsilon|G|$ neighbours, or there are two
disjoint sets $A,B\subseteq V(G)$ with $|A|,|B|\ge \epsilon|G|$ such that no edge joins $A$ and $B$. It follows 
that for every graph $H$, there exists $c>0$ such that for every graph $G$, if no induced subgraph of $G$ or its 
complement is a subdivision of $H$, then $G$ has a clique or stable set of cardinality at least $|G|^c$.
This is related to the Erd\H{o}s-Hajnal conjecture.
\end{abstract}

\section{Introduction}

For a graph $G$, we write $\omega(G), \alpha(G)$ for the cardinalities of the largest clique and largest stable set in $G$
respectively.  The number of vertices of $G$ is denoted by $|G|$, and $\overline{G}$ denotes the complement graph of $G$.
If $v\in V(G)$, $N(v)$ denotes the set of neighbours of $v$. Subsets $A,B$ of $V(G)$ are 
{\em complete} if $A\cap B=\emptyset$ and every vertex of $A$ is adjacent to every vertex of $B$, and 
{\em anticomplete} if $A\cap B=\emptyset$ and no vertex in $A$ has a neighbour in $B$. A pair $(A,B)$ of subsets of 
$V(G)$ is {\em pure} if 
$A$ is either complete or anticomplete to $B$.
For graphs $G,H$, we say $G$ is {\em $H$-free} if no induced subgraph of $G$ is isomorphic to $H$.
(All graphs in this paper are finite and have no loops 
or parallel edges.) An {\em ideal} of graphs is a class of graphs closed under isomorphism and under taking induced subgraphs; and an 
ideal is {\em proper} if it is not the class of all graphs.

It is well-known from Ramsey theory \cite{ES} that every graph $G$ contains a clique or stable set of size at least 
$\frac12\log |G|$.  On the other hand, there are graphs $G$ with no clique or stable set of size more than $2\log |G|$ \cite{E0} 
(in fact, most graphs have this property).  The celebrated Erd\H{o}s-Hajnal conjecture asserts that $H$-free graphs have {\em much} larger cliques or stable sets.  
Let us say that an ideal $\mathcal I$  has the {\em Erd\H{o}s-Hajnal property} if
there is some $\epsilon>0$ such that every graph $G\in\mathcal I$ has a clique or stable set of size at least $|G|^\epsilon$.  
The Erd\H os-Hajnal conjecture \cite{EH0,EH} is the following:

\begin{thm}\label{EHconj0}
{\bf Conjecture: }For every graph $H$, the ideal of $H$-free graphs has the Erd\H os-Hajnal property.
\end{thm}


A related but stronger property for an ideal is that every graph in the ideal contains a pure pair of linear-sized sets.
More formally, 
let us say that an ideal $\mathcal I$  has the {\em strong Erd\H os-Hajnal property} if
there is some $\epsilon>0$ such that every graph $G\in\mathcal I$ with at least two vertices 
contains a pure pair of sets that both have size at least $\epsilon |G|$.
It is easy to show that if an ideal has the strong Erd\H{o}s-Hajnal property then it has the Erd\H{o}s-Hajnal property 
(see \cite{APPRS,fp}; or section \ref{rt} below).
But the reverse implication does not hold.  In fact, if the ideal of $H$-free graphs has the strong Erd\H{o}s-Hajnal property then both $H$ and $\overline H$ are forests 
(to show that $H$ must be a forest, suppose that $H$ contains a cycle of length $k$, take a 
random graph $G\in{\mathcal G}(n,p)$, where $p$ is chosen so that $np\to \infty$ and $(np)^{k}=o(n)$, and delete one vertex from each cycle of length $k$; for $\overline H$, take complements).
Thus the ideal of $H$-free graphs does not have the strong Erd\H{o}s-Hajnal property for any graph $H$ with more than four vertices.
In this paper, we are interested in which ideals do have the strong Erd\H{o}s-Hajnal property.

An ideal is characterized by the minimal induced subgraphs that it does not contain.  If an ideal is defined by a finite number of excluded induced subgraphs, then the random graph construction shows that one of them must be a forest and one of them must be the complement of a forest.  An important result of this type 
 is due to Bousquet,  Lagoutte and Thomass\'e~\cite{lagoutte}.  Improving on earlier work of Chudnovsky and Zwols~\cite{cz}  and Chudnovsky and Seymour~\cite{cs0},
they showed that for every path $P$, the ideal of graphs with no induced $P$ or $\overline P$ has the strong Erd\H{o}s-Hajnal property.  
More recently, Liebenau, Pilipczuk, Seymour and Spirkl~\cite{caterpillars} (improving on an earlier result of Choromanski, Falik, Patel and Pilipczuk~\cite{hooks}) showed that
if $T$ is a subdivision of a caterpillar then the ideal of graphs with no induced $T$ or $\overline T$ has the strong Erd\H{o}s-Hajnal property.   They further conjectured that
the same statement holds for any forest $T$.   This conjecture is proved by the current authors in~\cite{treesanti}, 
completing the classification for ideals defined by a finite number of excluded induced subgraphs.

What about ideals that are not defined by a finite number of excluded induced subgraphs?  A breakthrough result in this direction was proved by
Bonamy, Bousquet and Thomass\'e~\cite{bonamy}, who showed that for every $k$ the ideal of graphs $G$ such that neither $G$ nor $\overline G$ contains an induced cycle
of length at least $k$ has the strong Erd\H{o}s-Hajnal property.  In other words, we exclude induced subdivisions of the cycle $C_k$ from both $G$ and $\overline G$.  In this paper, we prove a very substantial extension of this result.

\begin{thm}\label{result0}
For every graph $H$, the ideal of graphs $G$ such that neither $G$ nor $\overline G$ contains an induced subdivision of $H$
has the strong Erd\H{o}s-Hajnal property.
\end{thm}
If instead we take the ideal of all graphs $G$ that do not contain an induced subdivision of $H$, 
then in general such ideals need not have the strong
strong Erd\H{o}s-Hajnal property. For instance, the ideal of all graphs that do not contain an induced subdivision 
of a cycle of length six does not have the strong Erd\H{o}s-Hajnal property, because it includes the ideal of complements of all
triangle-free graphs.

As an immediate corollary of \ref{result0} we obtain the following.

\begin{thm}\label{EHapp0}
For every graph $H$, there exists $c>0$ such that for every graph $G$, one of the following holds:
\begin{itemize}
 \item $G$ or its complement contains an induced subdivision of $H$;
 \item $G$ contains a clique or stable set of size at least $|G|^c$.
\end{itemize}
\end{thm}

We can say a little more than \ref{result0}.  
If an ideal satisfies the strong Erd\H{o}s-Hajnal property then we know that every graph in the ideal has a pair of large sets that are either complete or anticomplete,
but we may not be able to choose which (for instance, consider the ideal consisting of all vertex-disjoint unions of cliques).  However, a theorem of R\"odl \cite{rodl} (discussed in 
section \ref{rt}) allows us to
assume that our graph is either quite sparse or quite dense.  We can then deduce \ref{result0} from the following significantly stronger ``one-sided'' result.

\begin{thm}\label{EHapp}
For every graph $H$, there exists $c>0$ such that every graph $G$ with at least two vertices and at most $c|G|^2$ edges satisfies one of the following:
\begin{itemize}
 \item $G$ contains an induced subdivision of $H$;
 \item there are two anticomplete subsets of $V(G)$, both of size at least $c|G|$.
\end{itemize}
\end{thm}

In fact, we will prove an even more general result: we will show it is enough to consider induced subdivisions where the edges in a 
specified path are not subdivided; and we will prove a version of the result (stated as \ref{massthm}) that works when the graph is 
weighted.  We introduce the necessary definitions and state results formally in the next section.  We discuss R\"odl's theorem and its 
application in section \ref{rt}, and then give the proof of \ref{massthm} over the next four sections.  An important feature of the 
proof is to divide the problem into two cases: in one case, we may assume that all small balls have small mass; in the other, we may 
assume that a significant mass is always concentrated in a small ball.  After some initial work, these cases are handled separately in 
sections \ref{sec:bigballs} and \ref{moreballs}.  We conclude, in the final section, with some applications, and a discussion of the 
relationship between the Erd\H{o}s-Hajnal conjecture and questions about $\chi$-boundedness.

\section{Statement of results}
 
Every proper ideal is contained in the ideal of $H$-free graphs for some $H$.
Thus \ref{EHconj0} can be reformulated as:
\begin{thm}\label{EHconj}
{\bf Conjecture: }For every proper ideal $\mathcal{I}$, there exists $c>0$ such that every graph $G\in \mathcal{I}$ satisfies
$\omega(G)\alpha(G)\ge |G|^c.$
\end{thm}

For $\epsilon>0$, let us say a graph $G$ is {\em $\epsilon$-coherent} if 
\begin{itemize}
\item $|G|\ge 2$;
\item $|N(v)|< \epsilon|G|$ for each $v\in V(G)$; and
\item $\min(|A|,|B|)< \epsilon|G|$, for every two anticomplete sets $A,B\subseteq V(G)$.
\end{itemize}
As we explain in section \ref{rt}, \ref{EHconj} is equivalent to the following:

\begin{thm}\label{sparseconj}
{\bf Conjecture: }For every proper ideal $\mathcal{I}$  there exist $\epsilon>0$ and $c>0$ such that every $\epsilon$-coherent graph 
$G\in \mathcal{I}$ satisfies $\omega(G)\alpha(G)\ge |G|^c.$
\end{thm}

Let us say an ideal is {\em incoherent} if for some $\epsilon>0$, no member of $\mathcal{I}$ is
$\epsilon$-coherent; and {\em coherent} if there is no such $\epsilon$. 

Let $H$ be a graph and let $P$ be a subgraph of $H$. Let $J$ be a graph obtained from $H$ by subdividing at least once
every edge of $H$ not in $E(P)$, and not subdividing the edges in $E(P)$. We call such a graph $J$ (and graphs isomorphic to it) a 
{\em $P$-filleting} of $H$. Our main result states:
\begin{thm}\label{mainthm}
Let $H$ be a graph and let $P$ be a path of $H$. Then every coherent ideal contains a 
$P$-filleting of $H$. 
\end{thm}

Here are some consequences of \ref{mainthm}.
\begin{itemize}
\item By setting $H=K_t$ and $|P|=1$, it follows that
the ideal of graphs with no
induced subgraph a subdivision of $K_t$
is incoherent.
\item Let $P$ be a path of $H$, let $k\ge 1$ be an integer, and let $H_k$ be obtained from $H$ by subdividing every edge not in
$E(P)$ exactly $k$ times. Since every coherent ideal contains a $P$-filleting of $H_k$, it follows that every coherent ideal
contains a $P$-filleting of $H$ where every edge not in $P$ is
subdivided at least $k$ times.
\item If $T$ is a tree such that some path $P$ of $T$ contains all vertices of degree at least three (such a tree
is called a {\em caterpillar subdivision}),
it follows that every coherent ideal contains $T$
(this is the main theorem of~\cite{caterpillars}). To see this, observe that since $T$ is a
caterpillar subdivision, every $P$-filleting of $T$ contains a copy of $T$ as an induced subgraph.
\end{itemize}

It might be possible to strengthen \ref{mainthm}, to the following:
\begin{thm}\label{treeconj}
{\bf Conjecture: }Let $H$ be a graph and let $P$ be a forest of $H$. Then every coherent ideal contains a
$P$-filleting of $H$.
\end{thm}
If so, this would be best possible, in the sense that the same conclusion does not hold if $P$ contains a cycle of $H$ (because
there are coherent ideals in which every graph has girth at least any fixed integer). It would imply the main theorem
of~\cite{treesanti} in the same way that \ref{mainthm} implies the main theorem of~\cite{caterpillars}.

For every subset $X\subseteq V(G)$, let $\mu(X)$ be some real number, satisfying
\begin{itemize}
\item $\mu(\emptyset)=0$ and $\mu(V(G))=1$, and $\mu(X)\le \mu(Y)$ for all $X,Y$ with $X\subseteq Y$; and
\item $\mu(X\cup Y)\le \mu(X)+\mu(Y)$ for all disjoint sets $X,Y$.
\end{itemize}
We call such a function $\mu$ a {\em mass} on $G$, and we call the pair $(G,\mu)$ a {\em massed graph}. 
For instance, we could take $\mu(X)=|X|/|G|$, or
$\mu(X) = \chi(G[X])/\chi(G)$, where $\chi$ denotes chromatic number. 
(It is sometimes convenient to speak of the ``mass'' of a set $X$, meaning $\mu(X)$.)
For $\epsilon>0$ let us say a massed graph $(G,\mu)$ is {\em $\epsilon$-coherent} if 
\begin{itemize}
\item $\mu(\{v\})< \epsilon$ for every vertex $v$; 
\item $\mu(N(v))< \epsilon$ for every vertex $v$; and
\item $\min(\mu(A),\mu(B))<\epsilon$ for every two anticomplete sets of vertices $A,B$.
\end{itemize}

In~\cite{caterpillars} it was found that
the main theorem of that paper could be extended to massed graphs: that for every caterpillar subdivision $T$,
there exists $\epsilon>0$ such that for every $\epsilon$-coherent massed graph $(G,\mu)$, some induced subgraph 
of $G$ is isomorphic to $T$.
We will show that \ref{mainthm} admits the same extension to masses. We will prove:

\begin{thm}\label{massthm}
For every graph $H$ and path $P$ of $H$, there exists $\epsilon>0$ such that
for every $\epsilon$-coherent massed graph $(G,\mu)$, some induced subgraph 
of $G$ is a $P$-filleting of $H$.
\end{thm}
\noindent{\bf Proof of \ref{mainthm}, assuming \ref{massthm}.\ \ }
Let $H$ be a graph and let $P$ be a path of $H$. Let $\mathcal{I}$ be a coherent ideal;
we must show that $\mathcal{I}$ contains a $P$-filleting of $H$. Choose $\epsilon>0$ satisfying \ref{massthm}.
By reducing $\epsilon$, we may assume that $\epsilon<1/2$.
Since $\mathcal{I}$ is coherent, some $G\in \mathcal{I}$ is $\epsilon$-coherent. Define $\mu(X)=|X|/|G|$ for
every $X\subseteq V(G)$. Then $(G,\mu)$ is a massed graph. We claim it is $\epsilon$-coherent. To see this we must
check the three conditions in the definition of $\epsilon$-coherence for massed graphs. The second and third follow
immediately from the second and third conditions in the definition of $\epsilon$-coherence for graphs, but the first
is not so clear. For the first we must show that 
$\mu(\{v\})< \epsilon$ for every vertex $v$; that is, $\epsilon|G|>1$. To see this, there are two cases. If all vertices in $G$
are pairwise adjacent, then since 
$|N(v)|< \epsilon|G|$ for each $v\in V(G)$, it follows that $|G|-1<\epsilon|G|$, and so $(1-\epsilon)|G|<1$, which is impossible
since $|G|\ge 2$ and $\epsilon\le 1/2$. If some two vertices $u,v$ of $G$ are nonadjacent, then $\{u\},\{v\}$ are anticomplete
sets, and so $1<\epsilon|G|$ since $G$ is $\epsilon$-coherent. This proves that  $(G,\mu)$ is $\epsilon$-coherent.
Consequently, from \ref{massthm},  some induced subgraph
of $G$ is a $P$-filleting of $H$. This proves \ref{mainthm}.~\bbox

\section{R\"odl's theorem}\label{rt}

Before we go on, let us prove the equivalence of \ref{EHconj} and \ref{sparseconj}. 
This is routine, and no doubt well-known to those familiar with the field, but we give the proof anyway.
Certainly \ref{EHconj} 
implies \ref{sparseconj}, and for the converse we use an invaluable tool due to 
R\"odl~\cite{rodl}, the following. 
\begin{thm}\label{rodl}
For every graph $H$ and all $\epsilon>0$ there exists $\delta>0$ such that for every $H$-free graph $G$,
there exists $X\subseteq V(G)$ with $|X|\ge \delta|G|$ such that in one of $G[X]$, $\overline{G}[X]$, 
every vertex in $X$ has degree less than $\epsilon|X|$.
\end{thm}

\ref{rodl} implies:
\begin{thm}\label{fullstrong}
If $\mathcal{I}$ is an ideal such that $\mathcal{I}$ and the ideal of complements of members of $\mathcal{I}$ are both
incoherent, then
there exists $c>0$ such that for all $G\in \mathcal{I}$ with $|G|>1$, there is a pure pair $(A,B)$ in $G$ with $|A|,|B|\ge c|G|$.
\end{thm}
\Proof
Let $\mathcal{I}_2$ be the ideal of complements of members of $\mathcal{I}$, and choose $\epsilon>0$ such that no member of $\mathcal{I}\cup \mathcal{I}_2$
is $\epsilon$-coherent. Choose a graph $H$ not in $\mathcal{I}$; and choose $\delta>0$ to satisfy \ref{rodl}, with
$H,\epsilon$ as given. Let $c=\delta\epsilon$. Now let $G\in \mathcal{I}$ with $|G|>1$.
We claim there is a pure pair $(A,B)$ in $G$ with $|A|,|B|\ge c|G|$. If $c|G|\le 1$, then we may take $A,B$ to be
disjoint singleton sets, and this is possible since $|G|>1$. Thus we may assume that $c|G|>1$.
Since $G$ does not contain $H$, by the choice of $\delta$
there exists $X\subseteq V(G)$ with $|X|\ge \delta|G|$ such that in one of $G[X]$, $\overline{G}[X]$, say $G'$,
every vertex in $X$ has degree less than $\epsilon|X|$. Thus $|X|\ge \delta|G|\ge c|G|>1$.
Since $G'\in \mathcal{I}\cup \mathcal{I}_2$, it follows that
$G'$ is not $\epsilon$-coherent; and so there exist $A,B\subseteq X$, anticomplete, with $|A|,|B|\ge \epsilon|X|\ge c|G|$,
as claimed. This proves \ref{fullstrong}.~\bbox

Another consequence of \ref{rodl} is:
\begin{thm}\label{rodlapp}
Let $\mathcal{I}$ be a proper ideal, let $\epsilon,c_0>0$, and suppose that for every $G\in \mathcal{I}$, if one of $G,\overline{G}$ is $\epsilon$-coherent
then 
$\omega(G)\alpha(G)\ge |G|^{c_0}$. Then $\mathcal{I}$ satisfies \ref{EHconj}.
\end{thm}
\Proof Let
$\mathcal{I}_2$ be the ideal of complements of members of $\mathcal{I}$.
Then every $\epsilon$-coherent graph
$G\in \mathcal{I}\cup \mathcal{I}_2$ satisfies $\omega(G)\alpha(G)\ge |G|^{c_0}.$
Choose $H\notin \mathcal{I}$;  
choose $\delta$ such that \ref{rodl} holds; and choose $c$ with $0<c\le c_0/2$
such that $\delta^{2c}\ge 1/2$ and $(\epsilon\delta)^c\ge 1/2$.
We prove by induction on $|G|$ that for
every graph $G\in \mathcal{I}\cup \mathcal{I}_2$, we have $\omega(G)\alpha(G)\ge |G|^c$.
If $|G|\le 1$ the claim is trivial, and if $2\le |G|\le \delta^{-2}$ then the claim holds, since
$$\omega(G)\alpha(G)\ge 2\ge \delta^{-2c}\ge |G|^c.$$
Thus we may assume that $|G|>\delta^{-2}$.
By \ref{rodl}, since one of $G,\overline{G}$ is $H$-free,
there exists $X\subseteq V(G)$ with $|X|\ge \delta|G|$ such that in one of $G[X]$, $\overline{G}[X]$,
every vertex in $X$ has degree at most $\epsilon|X|$. By replacing $G$ by $\overline{G}$ if necessary, we may assume that
every vertex in $X$ has degree at most $\epsilon|X|$ in $G[X]$.

Now $(\delta|G|)^{c_0}\ge (\delta|G|)^{2c}\ge |G|^c$, since $c_0\ge 2c$ and $|G|>\delta^{-2}$.
Thus if $G[X]$ is $\epsilon$-coherent, then
$$\omega(G[X])\alpha(G[X])\ge |X|^{c_0}\ge (\delta|G|)^{c_0}\ge |G|^c$$
as required. If $G[X]$ is not $\epsilon$-coherent, there exist two anticomplete subsets $A,B$ of $X$ such that
$|A|,|B|\ge \epsilon |X|$.
By the inductive hypothesis, $\omega(G[A])\alpha(G[A])\ge |A|^c,$ and the same for $B$,
and since $\alpha(G)\ge \alpha(G[A])+\alpha(G[B])$ and $\omega(G)\ge \omega(G[A]), \omega(G[B])$, it follows that
$$\omega(G)\alpha(G)\ge \omega(G[A])\alpha(G[A])+\omega(G[B])\alpha(G[B])\ge |A|^c+|B|^c\ge 2(\epsilon|X|)^c\ge 2(\epsilon\delta|G|)^c \ge |G|^c.$$
This proves \ref{rodlapp}.~\bbox

Let us show that if \ref{sparseconj} holds (for all proper ideals), then so does \ref{EHconj}. 
Let $\mathcal{I}$ be a proper ideal, and let
$\mathcal{I}_2$ be the ideal of complements of members of $\mathcal{I}$.
By applying \ref{sparseconj} to $\mathcal{I}$ and to $\mathcal{I}_2$, there exist $\epsilon,c_0$
such that every $\epsilon$-coherent graph
$G\in \mathcal{I}$ satisfies $\omega(G)\alpha(G)\ge |G|^{c_0},$ and
the same for $\mathcal{I}_2$. 
But then the result follows from \ref{rodlapp}.

Another useful consequence of \ref{rodlapp} is the following:
\begin{thm}\label{complements}
If $\mathcal{I}$ is an ideal such that $\mathcal{I}$ and the ideal of complements of members of $\mathcal{I}$ are both incoherent, then
$\mathcal{I}$ satisfies \ref{EHconj}.
\end{thm}
\Proof Let $\mathcal{I}_2$ be the ideal of complements of members of $\mathcal{I}$, and choose $\epsilon>0$ such that no member of $\mathcal{I}\cup \mathcal{I}_2$
is $\epsilon$-coherent. Then the result follows from \ref{rodlapp}.~\bbox

Next let us deduce the claims of section 1. It is possible to prove 
\ref{EHapp} $\Rightarrow$ \ref{result0} $\Rightarrow$ \ref{EHapp0}, as we said in the introduction, but 
it is more convenient to derive them all directly from results proved in this section.

\bigskip

\noindent{\bf Proof of \ref{result0}, assuming \ref{mainthm}.\ \ }
Let $H$ be a graph, and let $\mathcal{I}$ be the ideal of all graphs $G$ such that neither $G$ nor $\overline G$ contains
an induced subdivision of $H$. Thus $\mathcal{I}=\mathcal{I}_2$,
where $\mathcal{I}_2$ is the ideal of complements of members of $\mathcal{I}$.
By \ref{mainthm}, $\mathcal{I}$ is incoherent. By \ref{fullstrong}, 
there exists $c>0$ such that for all $G\in \mathcal{I}$ with $|G|>1$, there is a pure pair $(A,B)$ in $G$ with $|A|,|B|\ge c|G|$.
Hence $\mathcal{I}$ has the strong Erd\H{o}s-Hajnal property. This proves \ref{result0}.~\bbox

\noindent{\bf Proof of \ref{EHapp0}, assuming \ref{mainthm}.\ \ }
Let $H$ be a graph, and let $\mathcal{I}$ be the ideal of all graphs $G$ such that neither $G$ nor $\overline G$ contains
an induced subdivision of $H$. As before, $\mathcal{I}=\mathcal{I}_2$,
where $\mathcal{I}_2$ is the ideal of complements of members of $\mathcal{I}$, and $\mathcal{I}$ is incoherent.
By \ref{complements}, $\mathcal{I}$ satisfies \ref{EHconj}.~\bbox

\noindent{\bf Proof of \ref{EHapp}, assuming \ref{mainthm}.\ \ }
Let $H$ be a graph, and let $\mathcal{I}$ be the ideal of all graphs that contain no induced subdivision of $H$.
Let $P$ be a one-vertex path of $H$. By \ref{mainthm}, every coherent ideal contains a
$P$-filleting of $H$, and so $\mathcal{I}$ is incoherent. Choose $\epsilon>0$ such that 
no member of $\mathcal{I}$ is $\epsilon$-coherent, and let $c=\epsilon/9$. We claim that $c$ satisfies \ref{EHapp}.
Let $G$ be a graph with $|G|>1$ and $|E(G)|\le c|G|^2$ that contains no induced subdivision of $H$. 
Let $Y$ be the set of vertices of $G$ with degree at least
$\epsilon |G|/2$. Then $|Y|\epsilon|G|/2\le 2|E(G)|\le 2c|G|^2$, and so $|Y|< |G|/2$. Let $X=V(G)\setminus Y$;
so $|X|> |G|/2$, and so $|X|\ge 2$ since $|G|\ge 2$.
Every vertex of $G[X]$ has degree in $G[X]$ less than $\epsilon|G|/2\le \epsilon|X|$. Since $G[X]\in \mathcal{I}$,
$G[X]$ is not $\epsilon$-coherent, and so there exist anticomplete subsets $A,B$ of $X$ with $|A|,|B|\ge \epsilon|X|\ge c|G|$,
as required. This proves \ref{EHapp}.~\bbox

\section{Some preliminaries}

In order to prove \ref{massthm}, we might as well assume that $P$ is a Hamilton path of $H$. To see this, let 
$P$ have vertices $v_1\ll v_k$ in order and let the remaining vertices of $H$ be $v_{k+1}\ll v_n$. Add new vertices
$u_{k+1}\ll u_n$ to $H$, where each $u_i$ is adjacent to $v_{i-1}$ and $v_{i}$; let the new graph be $H'$
and let $P'$ be the path with vertices
$$v_1\ll v_k, u_{k+1}, v_{k+1}, u_{k+2}\ll u_n, v_n.$$
Then $P'$ is a Hamilton path of $H'$, and if the theorem holds for $(H', P')$ then it holds for $(H,P)$.
We will therefore assume that $P$ is a Hamilton
path of $H$.
In order to find a $P$-filleting of $H$ as an induced subgraph of $G$, we need to find an induced path 
$Q$ say of $G$, with the same number of vertices as $P$, such that certain pairs of vertices of $Q$ are joined
by induced paths in $G$, pairwise disjoint and disjoint from $Q$ (except for their ends), such that their union with $Q$ is induced
in $G$. (In particular, there must be no edges of $G$ between their interiors.) 

A {\em pairing} $\Pi$ in a graph $G$ is a set of pairwise disjoint subsets of $V(G)$, each of cardinality one or two; and let
$V(\Pi)$ be the union of the members of $\Pi$. If $X\subseteq V(G)$, a {\em pairing of $X$} means a pairing $\Pi$ with
$V(\Pi)= X$. A pairing $\Pi$ of $X$ is {\em feasible} in $G$ if for each $e\in \Pi$ with $|e|=2$
there is an induced path $P_e$ of $G$ joining the two members of $e$, and for each $e\in \Pi$ with $|e|=1$, $P_e$
is the one-vertex path with vertex set $e$, such that 
for all distinct $e,f\in \Pi$, the sets $V(P_e),V(P_f)$ are anticomplete.

An induced subgraph of $G$ isomorphic to $T$ is called a {\em copy} of $T$ in $G$.
Let us say a {\em caterpillar} is a tree
in which some path contains all vertices with degree more than one. Its
{\em leaves} are its vertices of degree one.  A {\em leaf-pairing} of $T$ means a pairing of the set of leaves of $T$.
A {\em caterpillar in $G$} means an induced subgraph of $G$
that is a caterpillar. 
Let $T$ be a caterpillar in $G$, and let $\Pi$ be a leaf-pairing of $T$.
Let $X$ be the set of all vertices in $V(G)\setminus V(T)$ with no neighbours in $V(T)\setminus V(\Pi)$.
The pairing $\Pi$
is {\em feasible} in $G$ {\em relative to $T$} if $\Pi$ is feasible in $G[X\cup V(\Pi)]$.
Thus,
another way to pose the problem of \ref{massthm} is to say that we are given a caterpillar $T$ with a leaf-pairing,
and we are searching for a copy $T'$ of $T$ in $G$
such that the corresponding leaf-pairing of $T'$ is feasible in $G$ relative to $T'$. 

Let $T$ be a caterpillar in $G$.
We say $T$ is {\em versatile} in $G$ if {\em every} leaf-pairing
of $T$ is feasible in $G$ relative to $T$.
In order to prove \ref{massthm} it therefore suffices to prove the following strengthening.

\begin{thm}\label{catversion}
For every caterpillar $T$, there exists $\epsilon>0$ such that 
for every $\epsilon$-coherent massed graph $(G,\mu)$, there is a versatile copy of $T$ in $G$.
\end{thm}

If $u$ is a vertex of a graph $H$, we denote by $N^r(u)$ the set of all vertices $v$ of $G$ such that the distance 
between $u,v$ is exactly $r$, and $N^r[u]$ the set of $v$ such that this distance is at most $r$.
For $r\ge 1$ an integer, and $\epsilon>0$, let us say a massed graph $(G,\mu)$ is 
{\em $(\epsilon,r)$-coherent} if 
\begin{itemize}
\item $\mu(N^r[v])< \epsilon$ for every vertex $v$; and
\item $\min(\mu(A),\mu(B))<\epsilon$ for every two anticomplete sets of vertices $A,B$.
\end{itemize}
(Thus, $(\epsilon,1)$-coherent $\Rightarrow$ $\epsilon$-coherent $\Rightarrow$ $(2\epsilon,1)$-coherent.)

The proof of \ref{catversion} breaks into two parts. Given a caterpillar $T$, we will first prove the statement of \ref{catversion}
for massed graphs $(G,\mu)$ that are $(\epsilon, r)$-coherent (for some appropriate value of $r$ depending on $T$
but not on $G$); and then we will use this to prove \ref{catversion} in general.
The first part is more difficult, and carried out in section~\ref{sec:bigballs}.

\section{Finding a caterpillar}

We need a result which is a modification of the main theorem (2.6) of \cite{caterpillars}. The main idea of its proof is
exactly that of~\cite{caterpillars}, but we need several minor changes, and it seemed best to prove the whole thing again.
First we need a lemma (also proved in \cite{caterpillars}). We say $X\subseteq V(G)$ is {\em connected} if $G[X]$
is connected.

\begin{thm}\label{components}
Let $(G,\mu)$ be an $\epsilon$-coherent massed graph, and let $Y\subseteq V(G)$ with $\mu(Y)\ge 3\epsilon$. 
Then there is a connected subset $X\subseteq Y$ with
$\mu(X)> \mu(Y)-\epsilon$.
\end{thm}
\Proof Let the vertex sets of the components of $G[Y]$ be $X_1\ll X_k$ say. Choose
$i\ge 1$ minimal such that $\mu(X_1\cup\cdots\cup X_i)\ge \epsilon$. Since the sets 
$X_1\cup\cdots\cup X_i$ and $X_{i+1}\cup \cdots\cup X_n$ are anticomplete, it follows that
$\mu(X_{i+1}\cup \cdots\cup X_n)<\epsilon$;
and from the minimality of $i$, $\mu(X_1\cup\cdots\cup X_{i-1})< \epsilon$. But
$$\mu(X_1\cup\cdots\cup X_{i-1}) + \mu(X_i) + \mu(X_{i+1}\cup \cdots\cup X_n)\ge \mu(Y)\ge 3\epsilon,$$
and so $\mu(X_i)\ge \epsilon$. Since the sets $X_i$ and $Y\setminus X_i$ are anticomplete, it follows that
$\mu(Y\setminus X_i)<\epsilon$, and so $\mu(X_i)>\mu(Y)-\epsilon$.
This proves \ref{components}.~\bbox

A {\em rooted caterpillar} is a tree $T$ with a distinguished vertex $h$, called its {\em head}, such that some path of $T$
with one end $h$ contains all the vertices with degree more than one. 
A rooted caterpillar $T$ with more than one vertex has a unique {\em predecessor} $T'$ (up to isomorphism), defined as follows.
Let $h$ be the head of $T$.
\begin{itemize}
\item 
If $h$ is adjacent to some leaf $u$ of $T$, let $T'$ be the rooted caterpillar
obtained from $T$ by deleting $u$, with the same head $h$.
\item 
If $h$ has no neighbours that are leaves, then $h$ is a leaf; let its neighbour be $u$, and let $T'$
be the rooted caterpillar obtained by deleting $v$, with head $u$.
\end{itemize}
Thus, every rooted caterpillar can be grown in canonical one-vertex steps from a one-vertex rooted caterpillar.
If $T$ is a rooted caterpillar with $n$ vertices, say, let $T_1\ll T_n$ be the rooted caterpillars such that
$T_n = T$, and $|T_1|= 1$, and $T_{i-1}$ is the predecessor of $T_i$ for $2\le i\le n$. We call $T_1\ll T_n$
the {\em ancestors} of $T$.

Let $\mathcal{Y}$ be a set of pairwise disjoint subsets of $V(G)$, where $G$ is a graph. Let $N$ be a graph,
and for each $v\in V(N)$ let $X_v\subseteq V(G)$. We say that the family $X_v\;(v\in V(N))$ is {\em $\mathcal{Y}$-spread}
if for each $v\in V(N)$ there exists $Y_v\in \mathcal{Y}$ such that the sets $Y_v\;(v\in V(N))$ are all different,
and $X_v\subseteq Y_v$ for each $v\in V(N)$.

If $A,B\subseteq V(G)$ are disjoint, we say $A$ {\em covers} $B$ if every vertex in $B$ has a neighbour in $A$.
Let $(G,\mu)$ be a massed graph. We say $X\subseteq V(G)$ is {\em $\delta$-dominant} if 
$\mu(X \cup \bigcup_{x \in X}N(x)) \geq \delta$.

The distance in $G$ between $u,v$ is called the
{\em $G$-distance between $u,v$}.
For $X\subseteq V(G)$ and $v\in X$, let us say $v$ is an {\em $r$-centre} of $X$ if every vertex in $X$
has $G[X]$-distance at most $r$ from $v$ (and consequently $X$ is connected).
Let us say a massed graph $(G,\mu)$ is {\em $(\delta,r)$-focussed} if
for every $Z\subseteq V(G)$ with $\mu(Z)\ge \delta$,
there is a vertex $v\in Z$ with $\mu(N^r_{G[Z]}[v])\ge \mu(Z)/2$. 

Let $N$ be the union of one or more rooted caterpillars with pairwise anticomplete vertex sets.
(Thus each component of $N$ has a head.) Let $H$ be the set of heads of the components of $N$.
Let $(G,\mu)$ be a massed graph.
A {\em $\delta$-realization} of $N$ in $G$ is an assignment of a subset $X_v\subseteq V(G)$ to each vertex $v\in V(N)$,
satisfying the following conditions:
\begin{itemize}
\item the sets $X_v\;(v\in V(N))$ are pairwise disjoint;
\item for every edge $uv$ of $N$, if $v$ lies on the path of $N$ 
between $u$ and the head of the component of $N$ containing $u$,
then $X_u$ covers $X_v$;
\item for all distinct $u,v\in V(N)$, if $u,v$ are nonadjacent in $N$ and not both in $H$ then $X_u, X_v$ are anticomplete; and
\item for each $v\in H$, $\mu(X_v) \ge \delta$, and for each $v\in V(N)\setminus H$, $X_v$ is connected and $\delta$-dominant.
\end{itemize}

\begin{thm}\label{findcat}
Let $T$ be a rooted caterpillar, let $\delta,\epsilon>0$,
let $(G,\mu)$ be an $\epsilon$-coherent massed graph, and let $\mathcal{Y}$ be a set of disjoint subsets
of $V(G)$ such that $|\mathcal{Y}|=2^{|T|}$ and $\mu(Y)\ge 2^{2^{|T|}}(\delta+\epsilon)$ for each $Y\in \mathcal{Y}$. 
Then there is 
a $\mathcal{Y}$-spread $\delta$-realization of $T$ in $G$.

If in addition $r\ge 0$ is an integer, $\epsilon\le  \delta/2$, and $(G,\mu)$ is $(\delta,r)$-focussed, then there is
a $\mathcal{Y}$-spread $\delta$-realization $(X_v\;:v\in V(T))$ of $T$ in $G$ such that $X_v$ has an $r$-centre
for each $v\in V(T)$
except the head.
\end{thm}
\Proof
There are two cases: in one (let us call this the ``focussed case'') we have the additional hypotheses that $r\ge 0$ is an integer, $\epsilon\le  \delta/2$, and 
$(G,\mu)$ is $(\delta,r)$-focussed; and in the other case (the ``unfocussed case'') we do not assume this.
Let $p = 2^{|T|}$ and for $0\le i\le p$ let $m_i = 2^{i}(\delta+\epsilon) - \epsilon$.
Thus $m_0 = \delta$, and
$m_{i+1} = 2m_i  + \epsilon$ for $0\le i< p$.

If $N$ is a disjoint union of rooted caterpillars, each isomorphic to an ancestor of $T$, we call $N$ a {\em nursery}, and 
we define 
$\phi(N)=\sum_C 2^{|C|}$, where the sum is taken over all components $C$ of $N$.
Let $N_p$ be the nursery with $p$ components, each an isolated vertex. Thus $\phi(N_p) = 2p$, and 
since $2^p(\delta+\epsilon)\ge m_p$, the members of $\mathcal{Y}$ form a $\mathcal{Y}$-spread $m_p$-realization of $N_p$ in $G$.
Choose $k\le p$ minimum such that there is a nursery $N_k$ with $k$ components and with $\phi(N_k)\ge 2p$, 
and there is a 
$\mathcal{Y}$-spread $m_k$-realization of $N_k$ in $G$.
Since $\phi(N_k)\ge 2p$ and each component of $N_k$ is an ancestor of $T$, it follows that $N_k$ has at least two components, 
and so $k\ge 2$. Suppose (for a contradiction) that 
each component of $N_k$ is isomorphic to an ancestor of $T$ different from $T$.

Let the components of $N_k$ be $H_1\ll H_k$, where $|H_1|\le \cdots \le |H_k|$,
and for $1\le i\le k$ let $h_i$ be the head of $H_i$. Now for $1\le i\le k$ 
there is an ancestor $S_i$ of $T$ such that $H_i$ is the predecessor
of $S_i$, since $H_i$ is not isomorphic to $T$. We recall that $S_i$ is obtained from $H_i$ by 
adding a new leaf adjacent to $h_i$, and either keeping the same head, or making the new vertex the new head.
Let $I$ be the set of all $i\in \{1\ll k\}$ such that $H_i, S_i$ 
have different heads.
If $I\ne \emptyset$, choose $i\in I$, maximum,
and otherwise let $i = 1$.

Let $(X_v\;:v\in V(N_k))$ be a $\mathcal{Y}$-spread $m_k$-realization of $N_k$ in $G$.
We will choose $Z\subseteq X_{h_i}$ with $\mu(Z)\ge \epsilon$, and an ordering $\{z_1\ll z_n\}$ of the elements of $Z$,
but we 
treat the foccussed and unfocussed cases differently. Suppose first we are in the unfocussed case.
Since $k\ge 2$ and hence $m_k\ge 3\epsilon$, \ref{components} implies that 
there exists $Z\subseteq X_{h_i}$ with $\mu(Z)>\mu(X_{h_i})-\epsilon\ge m_k-\epsilon$ such that $Z$ is connected.
Number the vertices of $Z$ as $z_1\ll z_n$ say, such that
$\{z_1\ll z_q\}$ is connected for $1\le q\le n$.

In the focussed case, since
$(G,\mu)$ is $(\delta,r)$-focussed and $\mu(X_{h_i})\ge m_k\ge \delta$, there exists
$Z\subseteq X_{h_i}$ with $\mu(Z)>\mu(X_{h_i})/2$ and with an $r$-centre.
Thus $\mu(Z)\ge \epsilon$
since $\mu(Z)\ge \mu(X_{h_i})/2\ge m_k/2 \ge \delta/2\ge \epsilon$, by hypothesis.
Let $z_1$ be an $r$-centre of $Z$, and choose the ordering $z_1\ll z_n$ of the vertices of $Z$
in increasing order of $G[Z]$-distance from $z_1$.

Since $k\ge 2$, there exists $j\ne i$ with $1\le j\le k$; and since $\mu(Z)\ge \epsilon$,
the set of vertices in $X_{h_j}$ with a neighbour in $Z$ has mass more than $\mu(X_{h_j})-\epsilon\ge m_{k-1}$.
Consequently we may choose $q$ with $0\le q\le n$, minimum such that for some $j\in \{1\ll k\}\setminus \{i\}$,
the set of vertices in $X_{h_j}$ with a neighbour in $\{z_1\ll z_q\}$ has mass at least $m_{k-1}$.
In particular, $\{z_1\ll z_q\}$ is $m_{k-1}$-dominant, and $q\ge 1$.

\begin{itemize}
\item If $j<i$, it follows that $i\in I$. Let $N_{k-1}$ be the graph obtained from $N_k$
by adding the edge $h_ih_j$, and deleting all vertices in $V(H_j)\setminus \{h_j\}$.
Let $H_i'$ be the component
of $N_{k-1}$ that contains the edge $h_ih_j$, and let us assign its head to be $h_j$.
Consequently $H_i'$ is isomorphic to $S_i$, and so 
$N_{k-1}$ is a nursery with $k-1$ components. Moreover,
$|H_i|\ge |H_j|$ (because $i>j$), and so
$\phi(N_{k-1})\ge \phi(N_k)$.
\item If $j>i$, it follows that $j\notin I$. Let $N_{k-1}$ be the graph obtained from $N_k$
by adding the edge $h_ih_j$, and deleting all vertices in $V(H_i)\setminus \{h_i\}$.
Let $H_j'$ be the component
of $N_{k-1}$ that contains the edge $h_ih_j$, and let us assign its head to be $h_j$.
Thus $H_j'$ is isomorphic to $S_j$, and again $N_{k-1}$ is a nursery with $k-1$ components and 
$\phi(N_{k-1})\ge \phi(N_k)$.
\end{itemize}
For each $v\in V(N_{k-1})$ define $X_v'$ as follows:
\begin{itemize}
\item if $v\ne \{h_1\ll h_k\}$ let $X_v'=X_v$;
\item let $X_{h_i}'=\{z_1\ll z_q\}$;
\item let $X_{h_j}'$ be the set of vertices in $X_{h_j}$ with a neighbour in $\{z_1\ll z_q\}$;
\item for $1\le \ell\le k$ with $\ell\ne i,j$, let $X_{h_{\ell}}'$ be the set of vertices in $X_{h_{\ell}}$ with no neighbour in
$\{z_1\ll z_q\}$.
\end{itemize}
We see that $X_{h_i}'$ covers $X_{h_j}'$, and has no edges to $X_{h_{\ell}}'$ for $1\le \ell\le k$ with $\ell\ne i,j$;
and $X_{h_i}'$ is connected and $m_{k-1}$-dominant; and in the focussed case, $X_{h_i}'$ has an $r$-centre.
Moreover, $\mu(X_{h_j}')\ge m_{k-1}$. Let $1\le \ell\le k$ with $\ell\ne i,j$; then, since $q\ge 1$ and
from the choice of $q$,
the mass of the set of vertices in $X_{h_{\ell}}$ with a neighbour in $\{z_1\ll z_{q-1}\}$ is less than $m_{k-1}$,
and hence $\mu(X_{h_{\ell}}\setminus X_{h_{\ell}}')<m_{k-1}+\epsilon$. Since $\mu(X_{h_{\ell}})\ge m_k$
and $m_k = 2m_{k-1} + \epsilon$, it follows that
$\mu(X_{h_{\ell}}')\ge m_{k-1}$. Thus $(X_v'\;:v\in V(N_{k-1}))$
is a $\mathcal{Y}$-spread $m_{k-1}$-realization of $N_{k-1}$ in $G$, contrary to the minimality of $k$ since
$\phi(N_{k-1})\ge \phi(N_k)\ge 2p$.

Consequently some component of $N_k$ is isomorphic to $T$; but then the theorem holds (since $m_k\ge \delta$).
This proves \ref{findcat}.~\bbox

\section{If no small ball has large mass}\label{sec:bigballs}

In this section we prove \ref{catversion} assuming that no ball with bounded radius has large mass.
We need first:

\begin{thm}\label{getcolumns}
Let $k\ge 0$ be an integer, and let $\epsilon,\kappa>0$
such that $\kappa+4\epsilon\le 1$ and $(k-1)\kappa\le 1$.
Let $(G,\mu)$ be an $\epsilon$-coherent massed graph.
Then there are $2k+1$ subsets $A_1\ll A_k, B_1\ll B_k, C$ of $V(G)$, 
pairwise disjoint, with the following properties:
\begin{itemize}
\item for $1\le i\le k$, $A_i$ is connected and covers $B_i$;
\item for $1\le i\le k$, $A_i, C$ are anticomplete;
\item for all distinct $i,j\in \{1\ll k\}$, $A_i$ is anticomplete to $A_j\cup B_j$;
\item $\mu(C)\ge 1-3k\epsilon$;  and
\item for $1\le i\le k$, the set of vertices in $C$ covered by $B_i$ has mass at least $\kappa-3k\epsilon$.
\end{itemize}
\end{thm}
\Proof
We proceed by induction on $k$; the result is trivial for $k=0$, taking $C=V(G)$, so we assume $k\ge 1$. 
Consequently we may assume that 
there are $2k-1$ subsets $A_1\ll A_{k-1}, B_1\ll B_{k-1}, C'$ of $V(G)$,
pairwise disjoint, with the following properties:
\begin{itemize}
\item for $1\le i\le k-1$, $A_i$ is connected and covers $B_i$;
\item for $1\le i\le k-1$, $A_i, C'$ are anticomplete;
\item for all distinct $i,j\in \{1\ll k-1\}$, $A_i$ is anticomplete to $A_j\cup B_j$; 
\item $\mu(C')\ge 1-3(k-1)\epsilon$; and
\item for $1\le i\le k-1$, the set of vertices in $C'$ covered by $B_i$ has mass at least $\kappa-3(k-1)\epsilon$.
\end{itemize}
Choose these subsets such that, in addition, $|B_1|+\cdots+|B_{k-1}|$ is minimum. For $1\le i\le k-1$, let $C_i$
be the set of vertices in $C'$ covered by $B_i$. Thus $\mu(C_i) \ge \kappa-3(k-1)\epsilon$, and 
from the minimality of $|B_1|+\cdots+|B_{k-1}|$, it follows that
$\mu(C_i)\le \kappa-(3k-4)\epsilon$.
Let $D = C'\setminus (C_1\cup\cdots\cup C_{k-1})$.
Thus 
$$\mu(D)\ge 1-3(k-1)\epsilon-(k-1)(\kappa-(3k-4)\epsilon)=(1-(k-1)\kappa) + (k-1)(3k-7)\epsilon.$$
We claim that $\mu(D)\ge 3\epsilon$. If $k=1$, the above implies that $\mu(D)=1$, and if $k=2$, the above implies that
$\mu(D)\ge 1-\kappa-\epsilon$; and so in either case $\mu(D)\ge 3\epsilon$, since $\kappa+4\epsilon\le 1$.
If $k\ge 3$, then $(k-1)(3k-7)\ge 3$ (indeed, $\ge 4$), and so the same displayed inequality implies that
$\mu(D)\ge 3\epsilon$ since $1-(k-1)\kappa\ge 0$. This proves the claim that $\mu(D)\ge 3\epsilon$.
Note that $D$ is anticomplete to $A_i\cup B_i$ for $1\le i<k$.

For $X\subseteq D$, let $B(X)$ denote the set of vertices in $C'\setminus X$ with a neighbour in $X$.
By \ref{components}, there exists a connected subset $X\subseteq D$ with $\mu(X)\ge \mu(D)-\epsilon$;
and hence there is a connected subset $X\subseteq D$ with $\mu(X\cup B(X))\ge \epsilon$.
Choose such a set $X$ minimal, and let $A_k=X$ and $B_k=B(X)$. 
Then $A_k$ is anticomplete to $A_i\cup B_i$ for $1\le i<k$, since $A_k=X\subseteq D$; and
$B_k$ is anticomplete to $A_i$ for $1\le i<k$, since $B_k=B(X)\subseteq C'$.

Choose $x\in A_k$ such that $A_k\setminus \{x\}$ is connected (or empty). 
Since $\mu(x)<\epsilon$
and $\mu(N(x))<\epsilon$, the minimality of $X$ implies that
$\mu(A_k\cup B_k)\le 3\epsilon$.
Let $C = C'\setminus (A_k\cup B_k)$.
Since $\mu(C')\ge 1-3(k-1)\epsilon$, it follows that $\mu(C)\ge 1-3k\epsilon$, and 
since the set of vertices
in $C$ with no neighbour in $A_k\cup B_k$ has mass less than $\epsilon$, it follows that the set $C_k$
say 
of vertices in $C$ with a neighbour in $A_k\cup B_k$ satisfies 
$\mu(C_k)\ge 1-(3k+1)\epsilon\ge \kappa-3k\epsilon$. 
Also $C_k$ is anticomplete to $A_k$, since $B_k = B(X)$, and so $B_k$ covers $C_k$.
For $1\le i<k$, $C_i\subseteq B_k\cup C$, and
$\mu(C_i\cap B_k)\le 3\epsilon$, and so 
$$\mu(C_i\cap C)\ge \mu(C_i)-3\epsilon\ge \kappa-3(k-1)\epsilon-3\epsilon = \kappa-3k\epsilon.$$
Since $B_i$ covers $C_i\cap C$,
this proves \ref{getcolumns}.~\bbox

This is used to prove the following. Let us say a {\em $k$-ladder} in a graph $G$
is a family of $3k$ subsets
$$A_1\ll A_k, B_1\ll B_k, C_1\ll C_k$$ 
of $V(G)$, pairwise disjoint and such that 
\begin{itemize}
\item for $1\le i\le k$, $A_i$ is connected and covers $B_i$, and $B_i$ covers $C_i$;
\item for $1\le i\le k$, $A_i, C_i$ are anticomplete; and
\item for all distinct $i,j\in \{1\ll k\}$, $A_i$ is anticomplete to $A_j\cup B_j\cup C_j$.
\end{itemize}
If in addition we have
\begin{itemize}
\item for $1\le i<j\le k$, $B_i$ is anticomplete to $C_j$
\end{itemize}
we say the ladder is {\em half-cleaned}. Let us say the {\em union} of the $k$-ladder is the triple $(A,B,C)$
where $A=\bigcup_{1\le i\le k}A_i$, $B=\bigcup_{1\le i\le k}B_i$,
and $C=\bigcup_{1\le i\le k}C_i$.

\begin{thm}\label{getcolumns2}
Let $\epsilon,\kappa>0$, and let $k\ge 0$ be an integer
such that $(k-1)k(\kappa+\epsilon)\le 1$ and $(k-1)(\kappa+\epsilon)+4\epsilon\le 1$.
Let $(G,\mu)$ be an $\epsilon$-coherent massed graph.
Then there is a half-cleaned $k$-ladder
$$A_1\ll A_k, B_1\ll B_k, C_1\ll C_k$$
in $G$ 
such that $\mu(C_i)\ge \kappa$ for $1\le i\le k$.
\end{thm}
\Proof Let $\kappa'=k(\kappa+\epsilon)$. Since $\kappa'+4\epsilon\le 1$ and $(k-1)\kappa'\le 1$, it follows from 
\ref{getcolumns} that
there are $2k+1$ subsets $A_1\ll A_k, B_1'\ll B_k', C$ of $V(G)$,
all disjoint, with the following properties:
\begin{itemize}
\item for $1\le i\le k$, $A_i$ is connected and covers $B_i'$;
\item for $1\le i\le k$, $A_i, C$ are anticomplete;
\item for all distinct $i,j\in \{1\ll k\}$, $A_i$ is anticomplete to $A_j\cup B_j'$; and
\item for $1\le i\le k$, the set of vertices in $C$ covered by $B_i'$ has mass at least $\kappa'$.
\end{itemize}
Inductively, suppose that $0\le j<k$ and
we have defined $B_1\ll B_j$ with $B_i\subseteq B_i'$ 
for $1\le i\le j$,
and we have defined disjoint subsets $C_1\ll C_j$ of $C$ such that 
for $1\le i\le j$, $B_i$ covers $C_i$ and is anticomplete
to $C\setminus (C_1\cup\cdots\cup C_i)$, with $\kappa\le \mu(C_i)\le \kappa+\epsilon$.  Thus 
$\mu(C_1\cup\cdots\cup C_j)\le (k-1)(\kappa+\epsilon)$, and since the set of vertices in $C$ covered by $B_{j+1}'$
has mass at least $\kappa'= k(\kappa+\epsilon)$, we may choose $B_{j+1}\subseteq B_{j+1}'$ minimal such that 
the set, $C_{j+1}$ say, of vertices
in $C\setminus (C_1\cup\cdots\cup C_j)$ covered by $B_{j+1}$ has mass at least $\kappa$. 
From the minimality of $B_{j+1}$, it follows that $\mu(C_{j+1})\le \kappa+\epsilon$. This completes the inductive definition
and so proves \ref{getcolumns2}.~\bbox

\begin{thm}\label{useladder}
Let $\epsilon>0$, and let $k\ge 0$ be an integer.
Let $(G,\mu)$ be an $(\epsilon,2)$-coherent massed graph, and let 
$$A_1\ll A_k, B_1\ll B_k, C_1\ll C_k$$
be a $k$-ladder in $G$
such that $\mu(C_i)\ge 3k\epsilon$ for $1\le i\le k$. Let $(A,B,C)$ be its union.
Suppose that $b_i\in B_i$ for $1\le i\le k$, such that $b_1\ll b_k$ are pairwise nonadjacent.
Then every pairing of $\{b_1\ll b_k\}$ is feasible in $G[A\cup B\cup C]$.
\end{thm}
\Proof For each $v\in B$, let $i(v)\in \{1\ll k\}$ such that $v\in B_{i(v)}$, and for $1\le i\le k$ let $D_i = A_i\cup B_i\cup C_i$.
Let $\Pi$ be a pairing of $\{b_1\ll b_k\}$, and let $\{s_1,t_1\}\ll \{s_n,t_n\}$ be the members of $\Pi$ with cardinality two.
For $1\le m\le n$ we will construct inductively a path $P_m$ between $s_m, t_m$ with the following properties:
\begin{itemize}
\item $V(P_m)\subseteq D_{i(s_m)}\cup D_{i(t_m)}$;
\item for $1\le \ell<m$, the sets $V(P_{\ell}), V(P_m)$ are anticomplete;
\item $V(P_m)$ is anticomplete to $\{b_1\ll b_k\}\setminus \{s_m,t_m\}$; and
\item at most two vertices of $V(P_m)$ belong to $C$, and at most two to $B\setminus \{s_m,t_m\}$.
\end{itemize}
Let $1\le m\le n$, and suppose we have constructed $P_1\ll P_{m-1}$; we construct $P_m$ as follows.
Let 
$$Z=\{b_1\ll b_k\}\cup ((V(P_1)\cup\cdots\cup V(P_{m-1}))\cap (B\cup C)).$$ 
Thus $|Z|\le k+4(m-1)\le 3k-4$, since $m\le n\le k/2$.
Let $X$ be the set of vertices in $C_{i(s_m)}$ that have $G$-distance at least three from every vertex of $Z$.
Since $(G,\mu)$ is $(\epsilon, 2)$-coherent, it follows that 
$\mu(X)\ge \mu(C_{i(s_m)})-(3k-4)\epsilon\ge \epsilon$. 
Let $Y$ be the set of vertices in $C_{i(t_m)}$ that have $G$-distance at least three from every vertex of $Z$;
then similarly $\mu(Y)\ge \epsilon$. Since $X\cap Y=\emptyset$ and
$(G,\mu)$ is $\epsilon$-coherent, there exist
$x\in X$ and $y\in Y$, adjacent. Since $B_{i(s_m)}$ covers
$C_{i(s_m)}$, there exists $x'\in B_{i(s_m)}$ adjacent to $x$; and since the distance between $x$ and $Z$ is at least three,
it follows that $x'$ has no neighbour in $Z$. Similarly there exists $y'\in B_{i(t_m)}$ adjacent to $y$. 
Since $A_{i(s_m)}$ is connected and covers $B_{i(s_m)}$, 
and similarly for $t_m$, it 
follows that the
subgraph of $G$ induced on $A_{i(s_m)} \cup A_{i(t_m)}\cup\{s_m,t_m,x,y,x',y'\}$ is connected. Choose an
induced path $P_m$ joining $s_m, t_m$ in this subgraph. Then $P_m$ satisfies the first and fourth bullets above.

We claim that for $1\le \ell<m$, $V(P_{\ell})$ and $V(P_m)$ are anticomplete.
Suppose not, and let $u\in V(P_{\ell})$ and $v\in V(P_m)$ be adjacent or equal. Now either $u$ is one of $s_{\ell},t_{\ell}$,
or $u$ belongs to one of $A_{i(s_{\ell})}, A_{i(t_{\ell})}$, or $u\in Z$; and either
$v$ is one of $s_{m},t_{m}$,
or $v$ belongs to one of $A_{i(s_{m})}, A_{i(t_{m})}$, or $v\in \{x,y,x',y'\}$. 
If $u\in A_{i(s_{\ell})}$, then all its neighbours in $A\cup B\cup C$ belong to $A_{i(s_{\ell})}\cup B_{i(s_{\ell})}$
from the definition of a $k$-ladder; and since $v$ is not in the latter set, it follows that $u\notin A_{i(s_{\ell})}$.
Similarly $u\notin A_{i(t_{\ell})}$, and $v\notin A_{i(s_{m})}\cup A_{i(t_{m})}$. Consequently $u,v\in B\cup C$.
Thus $u\in Z$, and so $v\notin \{x,y,x',y'\}$ from the choice of $x,y,x',y'$. Hence $v$ is one of $s_{m},t_{m}$. But
from the choice of $P_{\ell}$, $V(P_{\ell})$ is anticomplete to $\{b_1\ll b_k\}\setminus \{s_{\ell},t_{\ell}\}$,
a contradiction. Thus $P_m$ satisfies the second bullet.

For the third bullet, suppose that $u\in \{b_1\ll b_k\}\setminus \{s_m,t_m\}$ is adjacent to $v\in V(P_m)$. Since $u\in Z$,
it follows that $v\notin \{x,y,x',y'\}$; and as before $v\notin A_{i(s_{m})}\cup A_{i(t_{m})}$, and so $v$ is one of 
$s_m,t_m$, a contradiction since $b_1\ll b_k$ are pairwise nonadjacent. Thus $P_m$ satisfies the third bullet.

This completes the inductive definition, and hence proves \ref{useladder}.~\bbox

In turn, \ref{useladder} is used to prove the following, the main result of this section.

\begin{thm}\label{bigradius}
For every caterpillar $T$, there exist $\epsilon>0$ and an integer $r\ge 1$, such that
for every $(\epsilon,r)$-coherent massed graph $(G,\mu)$, there is a versatile copy of $T$ in $G$.
\end{thm}
\Proof We may assume that $|T|\ge 3$. Assign a head to $T$ to make it rooted, not one of the leaves.
Let $k=2^{|T|}$.
Let $r=5k$, and choose $\epsilon>0$ such that
$(k-1)k(2^k(3k+2)+1)\epsilon\le 1$.
Let $(G,\mu)$ be an $(\epsilon,r)$-coherent massed graph.
By \ref{getcolumns2}, there is a half-cleaned $k$-ladder
$$A_1\ll A_k, B_1\ll B_k, C_1\ll C_k$$
in $G$ such that $\mu(C_i)\ge 2^k(3k+2)\epsilon$ for $1\le i\le k$.
Then the unfocussed case of \ref{findcat} (taking $\delta = (3k+1)\epsilon$)  implies that there is
a $\{C_1\ll C_k\}$-spread $(3k+1)\epsilon$-realization $(X_v\;:v\in V(T))$ of $T$ in $G$.

Let $t_1\ll t_q$ be the vertices of $T$ that are not leaves, where $t_1$ is the head of $T$ and
$t_it_{i+1}$ is an edge of $T$ for $1\le i<q$.
Choose $x_1\in X_{t_1}$. For $2\le i\le q$ in turn, since $X_{t_{i}}$ covers $X_{t_{i-1}}$ by definition of a
realization, we may choose $x_{i}\in X_{t_{i}}$ adjacent to $x_{i-1}$. 
Since each $x_i$ belongs to one of $C_1\ll C_k$, it follows that $\{x_1\ll x_q\}$ is anticomplete to $A_1\cup\cdots\cup A_k$.
Also, since there are no edges between
$X_u, X_v$ for nonadjacent $u,v\in V(T)$, it follows that 
$x_1\ll x_q$ are the vertices in order of an induced path of $G$. 
For the same reason we have the following two statements:
\\
\\
(1) {\em For each leaf $v$ of $T$ with neighbour $t_j$ say, $x_j$ has a neighbour in $X_v$, and $X_v$ is anticomplete to 
$\{x_1\ll x_q\}\setminus \{x_j\}$.}
\\
\\
(2) {\em For all distinct leaves $u,v$ of $T$, $X_u$ is anticomplete to $X_v$.}

\bigskip

We recall that $C_1\ll C_k$ are pairwise disjoint, and $X_v\;(v\in V(T))$ is $\{C_1\ll C_k\}$-spread;
let $I$ be the set of $i\in \{1\ll k\}$ such that $X_v\subseteq C_i$ for some leaf $v$ of $T$. 
For each $i\in I$, let $v_i$ be the (unique) leaf of $T$ with $X_{v_i}\subseteq C_i$. Let $i\in I$, and let $j\in \{1\ll q\}$
such that $v_i$ is adjacent to $t_j$ in $T$. We define $x^i = x_j$. Thus there may be distinct values $i,i'\in I$
with $x^i=x^{i'}$.

Let $i\in I$. Since $X_{v_i}$
is $(2k+1)\epsilon$-dominant, and $(G,\mu)$ is $(\epsilon, r)$-coherent, there exists $u\in X_{v_i}$ with a neighbour
$v$ such that the $G$-distance between $v$ and $x_1$ is at least $r+1$. Consequently the $G$-distance between $u$ and $x_1$
is at least $r$. By (1), $X_{v_i}\cup \{x^i\}$ is connected, and so there is a path of $G[X_{v_i}\cup \{x^i\}]$
between $x^i$ and $u$; and hence there is a minimal path $P_i$ of $G[X_{v_i}\cup \{x^i\}]$ with one end $x^i$ and the other 
$u_i$ say, such that the $G$-distance between $x_1$ and $u_i$ is at least $q+4i$.
It follows that the $G$-distance between $x_1$ and $u_i$ is exactly $q+4i$.
Choose a vertex $b_i\in B_i$ adjacent to $u_i$. Let $c_i$ be the second vertex of $P_i$, that is, the vertex adjacent to 
$x^i$, and let $Q_i = P_i\setminus \{x^i\}$. Thus $Q_i$ is a path of $G[X_{v_i}]$.
The subgraph $T'$ of $G$ induced on $\{x_1\ll x_q\}\cup \{c_i:i\in I\}$ is a copy of $T$, by (1) and (2), and 
we will show that it is versatile.
\\
\\
(3) {\em For all distinct $i,j\in I$, the sets $V(Q_i)\cup \{b_i\}$ and $V(Q_j)\cup \{b_j\}$ are anticomplete.}
\\
\\
We may assume that $i<j$. 
Certainly $V(Q_i)$ and $V(Q_j)$ are anticomplete, by (2).
Also $b_i$ has no neighbour in
$V(Q_j)$ since the $k$-ladder is half-cleaned; so it remains to check that $b_j$ has no neighbour in $V(Q_i)\cup \{b_i\}$.
Let $v\in V(Q_i)$; then from the minimality of $Q_i$, the $G$-distance between
$v,x_1$ is at most $q+4i$. But the $G$-distance between $u_j$ and $x_1$ is $q+4j$, and
so the $G$-distance between $v,u_j$ is at least four. Consequently the
$G$-distance between $v,b_j$ is at least three, and in particular $b_j$ has no neighbour in $V(Q_i)$; and setting $v=u_i$,
since the $G$-distance between $v,b_j$ is at least three it follows that $b_i, b_j$ are nonadjacent. This proves (3).
\\
\\
(4) {\em For all $i\in I$, the sets $V(Q_i)\cup \{b_i\}$ and $\{x_1\ll x_q\}\setminus \{x^i\}$ are anticomplete.}
\\
\\
Let $j\in \{1\ll q\}$ and suppose that $x_j$ is adjacent to some $v\in V(Q_i)\cup \{b_i\}$ and $x_j\ne x^i$. 
Since the $G$-distance between $x_1,b_i$
is $q+4i$, and the $G$-distance between $x_1, x_j$ is at most $q-1$, it follows that $b_i, x_j$ are nonadjacent, and so
$v\in Q_i$, contrary to (1).

\bigskip

Let $Z$ be the set of vertices in $G$ with $G$-distance from $x_1$ at most $4k+q+4$, and for
$1\le i\le k$ let $C_i'=C_i\setminus Z$. Since $k\ge 2^{q+2}$, it follows that $r=5k\ge 4k+q+4$, and so $(G,\mu)$ is 
$(\epsilon,4k+q+4)$-coherent. Hence $\mu(Z)\le \epsilon$.
Thus for $1\le i\le k$, $\mu(C_i')\ge \mu(C_i)-\epsilon\ge 3k\epsilon$.
Let $Q =\bigcup_{i\in I}V(Q_i)\cup  \{b_i\}$.
Since every vertex in $Q$ has $G$-distance at most $4k+q+1$ from $x_1$, it follows that
every vertex in $C_i'$ has $G$-distance at least three from $Q$.
Let $B_i'$ be the set of vertices in $B_i\setminus \{b_i\}$ with no neighbours in $V(Q)$. Since every vertex in $C_i'$ has a neighbour in $B_i$
and has $G$-distance at least three from $Q$, it follows that $B_i'$ covers $C_i'$.
Hence the sets
$$A_i\;(i\in I), B_i'\cup \{b_i\}\;(i\in I), C_i'\;(i\in I)$$
form an $|I|$-ladder, with union $(A',B', C')$ say. 
Since $\mu(C_i')\ge 3k\epsilon$ and $r\ge 2$, 
every pairing of $\{b_i:i\in I\}$  is feasible in $G[A'\cup B'\cup C']$ by \ref{useladder}.

Since $A'\cup B'\cup C'\setminus \{b_i\;:i\in I\}$ is anticomplete to 
$$(Q\setminus \{b_i\;:i\in I\})\cup \{x_1\ll x_q\},$$ 
it follows from (3) and (4)
that every pairing of $\{c_i:i\in I\}$ is feasible in the subgraph induced on
$A'\cup B'\cup C'\cup \bigcup_{i\in I}V(Q_i)$. But $\{c_i:i\in I\}$ is the set of leaves of $T'$, and since 
$$A'\cup B'\cup C'\cup \bigcup_{i\in I}(V(Q_i)\setminus \{c_i\})$$
is anticomplete to $\{x_1\ll x_q\}$,
it follows that $T'$ is versatile.
This proves \ref{bigradius}.~\bbox

\section{The general proof}\label{moreballs}

Now we turn to the proof of \ref{catversion} in general. Fix the caterpillar $T$, and choose $\epsilon_r$ and $r$
to satisfy \ref{bigradius} with $\epsilon$ replaced by $\epsilon_r$. Now we will choose $\epsilon$ much smaller than $\epsilon_r$,
and try to prove that in every $\epsilon$-coherent massed graph $(G,\mu)$, some copy of $T$ is versatile.
We can therefore assume that for every $Z\subseteq V(G)$, there is no mass $\mu'$ on $G[Z]$ that is
$(\epsilon_r,r)$-coherent; and in particular (assuming $\mu(Z)>0$), the mass $\mu'$ on $G[Z]$ defined by $\mu'(X)=\mu(X)/\mu(Z)$
for each $X\subseteq Z$ is not $(\epsilon_r,r)$-coherent. Consequently, either there is a vertex
$v\in Z$ with $\mu(N^r_{G[Z]}[v])\ge \epsilon_0 \mu(Z)$, or there are two anticomplete sets $A,B\subseteq Z$
with $\mu(A),\mu(B)\ge \epsilon_r\mu(Z)$. The latter is only helpful if $\epsilon_r\mu(Z)\ge \epsilon$, but in that case
we can assume the latter never occurs. Thus, for every $Z\subseteq V(G)$ with $\mu(Z)\ge \epsilon/\epsilon_r$,
there is a vertex $v\in Z$ with $\mu(N^r_{G[Z]}[v])\ge \epsilon_r \mu(Z)$. Since $N^r_{G[Z]}[v]$
is anticomplete to $Z\setminus N^{r+1}_{G[Z]}[v]$, and $\mu'(N^r_{G[Z]}[v])\ge \epsilon_r$ 
(and because of the ``latter never occurs''
assumption above), it follows that
$\mu'(Z\setminus N^{r+1}_{G[Z]}[v])<\epsilon_r$, and so $\mu'(N^{r+1}_{G[Z]}[v])\ge 1-\epsilon_r$, that is,
$\mu(N^{r+1}_{G[Z]}[v])\ge (1-\epsilon_r)\mu(Z)$. Initially we could have chosen $\epsilon_r$ as small as we want, and 
in particular we may assume that $\epsilon_r\le 1/2$; and so $\mu(N^{r+1}_{G[Z]}[v])\ge \mu(Z)/2$. Because of this we will
be able to apply the focussed case of \ref{findcat}.

For $X\subseteq V(G)$ and $v\in V(G)$, we say $v$ {\em touches} $X$ if either $v\in X$
or $v$ has a neighbour in $X$; and otherwise $v$ is {\em anticomplete} to $X$.
We need the following.
\begin{thm}\label{joinleaves}
Let $t,r\ge 1$ be integers and $\epsilon>0$, and 
let $(G,\mu)$ be an $\epsilon$-coherent massed graph.
Let $T$ be a caterpillar in $G$ with $t$ vertices, and let $x_1\ll x_q$ be the vertices of $T$ with degree more than one.
For each leaf $v$ of $T$ let $x^v$ be its neighbour in $\{x_1\ll x_q\}$, and let $X_v\subseteq V(G)$, such that
\begin{itemize}
\item for each leaf $v$, $v\in X_v$, and $X_v\cap \{x_1\ll x_q\} = \emptyset$, and 
$X_v$ is anticomplete to $\{x_1\ll x_q\}\setminus \{x^v\}$;
\item for all distinct leaves $u,v$, $X_u$ is anticomplete to $X_v$;
\item for each leaf $v$, $x^v$ is an $r$-centre for $X_v\cup \{x^v\}$;
\item for each leaf $v$, $v$ is the unique neighbour of $x^v$ in $X_v$; and 
\item for each leaf $v$, $X_v\cup \{x^v\}$ is $(r+2)t^{t+1}\epsilon$-dominant.
\end{itemize}
Then $T$ is versatile.
\end{thm}
\Proof Define $\kappa_i=(r+2)t^{t-i+1}\epsilon$ for $0\le i\le t$. 
Let $\Pi$ be a pairing of the set of leaves $L$ of $T$. Let $L=\{v_1\ll v_{\ell}\}$, where
$\Pi$ consists of the sets $\{v_{2i-1},v_{2i}\}$ for $1\le i\le k$ for some $k\le \ell/2$, together with the singleton
sets $\{v_i\}$ for $2k+1\le i\le \ell$. Let $X^0_v=X_v\cup \{x^v\}$ for each $v\in L$. For $1\le i\le k$, we define
$X^i_v\;(v\in \{v_{2i+1}\ll v_{\ell}\})$ and $P_i$ inductively as follows. 
We assume $P_1\ll P_{i-1}$ and $X^{i-1}_v\;(v\in \{v_{2i-1}\ll v_{\ell}\})$ have been defined, such that
\begin{itemize}
\item for $1\le h\le i-1$, $P_h$ is an induced path between $v_{2h-1}, v_{2h}$, of length at most $2r+1$;
\item for $1\le h\le i-1$, $V(P_h)\setminus \{v_{2h-1}, v_{2h}\}$ is anticomplete to $V(T)\setminus \{v_{2h-1}, v_{2h}\}$;
\item for $1\le h<h' \le i-1$, $V(P_h)$ is anticomplete to $V(P_{h'})$;
\item for $1\le h\le i-1$ and $v\in \{v_{2i-1}\ll v_{\ell}\}$, $V(P_h)$ is anticomplete to $X^{i-1}_v$;
\item for $v\in \{v_{2i-1}\ll v_\ell\}$, $X^{i-1}_{v}$ is $\kappa_{i-1}$-dominant; and 
\item for $v\in \{v_{2i-1}\ll v_\ell\}$, $\{x^v,v\}\subseteq X^{i-1}_v$ and $x^v$ is an $r$-centre for $X^{i-1}_{v}$.
\end{itemize}
For each $w\in \{v_{2i+1}\ll v_{\ell}\}$, choose $X^i_w\subseteq X^{i-1}_w$,
minimal such that $X^i_w$ is $\kappa_i$-dominant and $x^w$ is a $r$-centre for $X^i_w$. By deleting a vertex in $X^i_w$
with maximum $G[X^i_w]$-distance from $x^w$, the minimality of $X^i_w$ implies that 
the set of vertices that touch $X^i_w$ has mass at most $\kappa_i+\epsilon$. Also the set of 
vertices that touch $\{x_1\ll x_q\}\cup \bigcup_{h<i}V(P_h)$ has mass at most $(q+(k-1)(2r+2))\epsilon$. 
Let $u=v_{2i-1}$ and $v=v_{2i}$.
Let $C$ be the set of all vertices that do not touch $\{x_1\ll x_q\}\cup \bigcup_{h<i}V(P_h)$ and do not touch
$X^i_w$ for $w\in \{v_{2i+1}\ll v_{\ell}\}$. 
Let $A,B\subseteq C$ be the sets of all vertices in $C$ that touch
$X^{i-1}_u$ and touch $X^{i-1}_v$ respectively.

Since $X^{i-1}_u$ is $\kappa_{i-1}$-dominant, it follows that 
$$\mu(A)\ge \kappa_{i-1}- (q+(k-1)(2r+2))\epsilon- \ell(\kappa_i+\epsilon).$$
The expression on the right side of this inequality is at least $\epsilon$,
since $q+\ell=t$ and $\ell\le t-1$ and $k\le \ell/2$ and $r\ge 1$ (we leave checking this to the reader);
and so $\mu(A)\ge \epsilon$.
The same holds for $\mu(B)$; and so $A,B$ are not anticomplete. Consequently there are vertices $a,b$,
adjacent or equal, such that $a$ touches $X^{i-1}_u$ and $b$ touches $X^{i-1}_v$, and $a,b$ are anticomplete to
$\{x_1\ll x_q\}\cup \bigcup_{h<i}V(P_h)$ and to $X^i_w$ for $w\in \{v_{2i+1}\ll v_{\ell}\}$.
Since $x^u$ is an $r$-centre for $X^{i-1}_u$, and $u$ is the unique neighbour of $x^u$ in $X^{i-1}_u$,
there is a path of $G[X^{i-1}_u\cup \{a\}]$ between 
$u$ and $a$ of length at most $r$, and the same for $v$; and therefore there is an induced path $P_i$ between 
$u, v$ of length at most $2r+1$, and all its vertices belong to $X^{i-1}_u\cup X^{i-1}_v\cup \{a,b\}$.
In particular, $V(P_i)$ is anticomplete to $V(P_h)$ for $h<i$, since $X^{i-1}_u\cup X^{i-1}_v\cup \{a,b\}$
is anticomplete to $V(P_h)$; and $V(P_i)$ is anticomplete to $X^{i}_w$ for $w\in \{v_{2i+1}\ll v_{\ell}\}$, 
since $X^{i-1}_u\cup X^{i-1}_v\cup \{a,b\}$ is complete to $X^{i}_w$. This completes the inductive definition.

But then the paths $P_1\ll P_k$ together with the singletons $\{c_i\}\;(2k+1\le i\le \ell)$ show that $\Pi$
is feasible relative to $T$. Consequently $T$ is versatile. This proves \ref{joinleaves}.~\bbox

This is used to prove:

\begin{thm}\label{smallradius}
Let $T$ be a caterpillar with $t$ vertices, let $r\ge 1$ be an integer, and let $\epsilon,\delta>0$ with $\epsilon\le \delta/2$, such that
$\delta\le 2^{-(t+2^t)}t^{-t}$ and 
$\epsilon\le 2^{-(t+2^t)}t^{-2t}(3r+5)^{-1}.$
Let $(G,\mu)$ be a $(\delta,r)$-focussed $\epsilon$-coherent massed graph. Then there is a versatile copy of $T$ in $G$.
\end{thm}
\Proof We may assume that $|T|\ge 3$, $T$ is rooted, and its head is an internal vertex.
Let $t=|T|$ and $k=2^t$. Let $\lambda = 1/k-\epsilon$.
\\
\\
(1) {\em There exist pairwise disjoint subsets $Y_1\ll Y_k$ of $V(G)$ with $\mu(Y_i)\ge \lambda$ for $1\le i\le k$.}
\\
\\
We define $Y_1\ll Y_k\subseteq V(G)$ inductively as follows. Let $1\le i<k$, and assume 
we have chosen $Y_1\ll Y_i\subseteq V(G)$, pairwise disjoint and each with 
$\lambda\le \mu(Y_i)\le \lambda+\epsilon$. Thus 
$$\mu(Y_1\cup\cdots\cup Y_i)\le (k-1)(\lambda+\epsilon),$$
and so 
$$\mu(V(G)\setminus (Y_1\cup\cdots\cup Y_i))\ge 1-(k-1)(\lambda+\epsilon)\ge \lambda.$$
Consequently we may choose $Y_{i+1}$ disjoint from $Y_1\cup\cdots\cup Y_i$ with $\mu(Y_i)\ge \lambda$. Choose $Y_{i+1}$
minimal with this property; then $\mu(Y_i)\le \lambda+\epsilon$. This completes the inductive definition of $Y_1\ll Y_k$
and so proves (1).

\bigskip
Let $\kappa_0=2^{-k}\lambda -\epsilon$, and 
for $1\le i\le t$ let $\kappa_i=\kappa_0 t^{-i}$. 
Since $G$ is $(\delta,r)$-focussed, it is also $(\kappa_0,r)$-focussed, since $\kappa_0\ge \delta$.
From the focussed case of \ref{findcat}, there is a $\{Y_1\ll Y_k\}$-spread $\kappa_0$-realization $(X_v\;:v\in V(T))$ of $T$ in $G$ 
such that $X_v$ has an $r$-centre for each $v\in V(T)$
except the head. Let the vertices of degree more than one in $T$ be $t_1\ll t_q$, where $t_1$ is the head and
$t_it_{i+1}$ is an edge of $T$ for $1\le i<q$. Choose $x_1\in X_{t_1}$, and inductively for $i = 2\ll q$, choose
$x_i\in X_{t_i}$ adjacent to $x_{i-1}$. As in the proof of \ref{bigradius}, $x_1\ll x_t$ are the vertices 
in order of an induced path of $G$. We need to arrange that for each leaf $v$ of $T$ adjacent to $t_i$ say, 
$x_i$ has a unique neighbour in $X_v$. 

Let $L=\{v_1\ll v_{\ell}\}$ be the set of leaves of $T$, and for each $v\in L$ let $x^v$ be the vertex 
$x_j$ such that $v$ is adjacent to $t_j$ in $T$.
Thus $x^v$ is the unique vertex in $\{x_1\ll x_q\}$ covered by $X_{v}$. Since $X_{v}$ has an $r$-centre, 
and $x^v$ has a 
neighbour in $X_{v}$, it follows that $x^v$ is a $(2r+1)$-centre of $X_{v}\cup \{x^v\}$.

Let $X^0_v=X_{v}\cup \{x^v\}$ for $v\in L$. For $0\le i\le \ell$
we will inductively define $X^i_v\;(v\in L)$ satisfying the following: for $0\le i\le \ell$,
\begin{itemize}
\item for each $v\in L$, $x^v\in X^{i}_v$ and $X^{i}_v\setminus\{x^v\}$ is anticomplete to $\{x_1\ll x_q\}\setminus \{x^v\}$;
\item for all distinct $u,v\in L$, $X^{i}_u\setminus \{x^u\}$ is anticomplete to $X^{i}_v\setminus \{x^v\}$;
\item for $1\le j\le i$, $x^{v_j}$ is a $(3r+2)$-centre of $X^{i}_{v_j}$, 
and $x^{v_j}$ has a unique neighbour in $X^{i}_{v_j}$; 
\item for $i< j\le \ell$, $x^{v_j}$ is a $(2r+1)$-centre of $X^{i}_{v_j}$; and
\item for each $v\in L$, $X^{i}_v$ is $\kappa_{i}$-dominant.
\end{itemize}
Suppose that $1\le i\le \ell$, and we have defined $X^{i-1}_v\;(v\in L)$ as above.
For each $u\in L\setminus \{v_i\}$, let $u = v_j$ say. Let $r_j = 3r+2$ if $j<i$, and $r_j=2r+1$ if $j>i$, so 
in either case $x^u$ is an $r_j$-centre of $X^{i-1}_u$. 
Choose $X^i_u\subseteq X^{i-1}_u$,
minimal such that $X^i_u$ is $\kappa_i$-dominant and $x^u$ is an $r_j$-centre for $X^i_u$. It follows 
from the minimality
of each $X^i_u$ that the set of vertices that touch $X^i_u$ has mass at most $\kappa_i+\epsilon$.

Let $v=v_i$, and let $Y''$ be the set of
all vertices that touch $X^{i-1}_v$; then $\mu(Y'')\ge \kappa_{i-1}$ since $X^{i-1}_v$ is $\kappa_{i-1}$-dominant.
Let $Y'$ be the set of vertices in $Y''$ that are anticomplete to 
$$\{x_1\ll x_q\}\cup \bigcup_{u\in L\setminus \{v\}}X^i_u.$$
Thus $\mu(Y')\ge \mu(Y'')-q\epsilon-(\ell-1)(\kappa_i+\epsilon)$. Since $\mu(Y'')\ge \kappa_{i-1}$ and
$\delta\le 2^{-(t+2^t)}t^{-t}$ by hypothesis, it follows that $\mu(Y')\ge \delta$ 
(we leave it to the reader to check this arithmetic).
Since $G$ is $(\delta,r)$-focussed,
there is a subset $Y\subseteq Y'$ with $\mu(Y)\ge \mu(Y')/2$, such that $Y$ has an $r$-centre $y$ say.
Since $x^v$ is a $(2r+1)$-centre for $X^{i-1}_v$, and $y$ touches this set, there is an induced path $P$ of 
$G[X^{i-1}_v\cup \{y\}]$ between
$x^v$ and $y$, of length at most $2r+2$. Let $X^i_v=Y\cup V(P)$.
Then $x^v$ is a $(3r+2)$-centre
for $X^i_v$. Since $x^v$ is anticomplete to $Y$ and has only one neighbour in $P$, it follows that
$x^v$ has only one neighbour in $Y\cup V(P)$. Moreover, $Y$ is $\kappa_i$-dominant since $\mu(Y)\ge \epsilon$.
This completes the inductive definition. 

For each $v\in L$, let $c_v$ be the unique neighbour of $x^v$ in $X^{\ell}_v$. The subgraph $T'$ induced on 
$\{x_1\ll x_q\}\cup \{c_v:\;v\in L\}$ is a copy of $T$. Since
$\epsilon\le 2^{-(t+2^t)}t^{-2t}(3r+5)^{-1}$ by hypothesis, it follows that
$\kappa_{\ell}\ge (3r+4)t^{t+1}\epsilon$, and so \ref{joinleaves} (with $r$ replaced by $3r+2$) implies
that $T'$ is versatile.
This proves \ref{smallradius}.~\bbox

Now let us put these pieces together, to prove \ref{catversion}, which we restate:
\begin{thm}\label{catversion2}
For every caterpillar $T$, there exists $\epsilon>0$, such that
for every $\epsilon$-coherent massed graph $(G,\mu)$, there is a versatile copy of $T$ in $G$.
\end{thm}
\Proof
By \ref{bigradius} there exist $\epsilon_r>0$ and an integer $r\ge 1$, such that
for every $(\epsilon_r,r)$-coherent massed graph $(G,\mu)$, there is a versatile copy of $T$ in $G$. We may assume 
$\epsilon_r\le t^{-t}(3r+5)^{-1}$. Choose $\epsilon$ such that 
$\epsilon\le 2^{-(t+2^t)}t^{-t}\epsilon_r$.
Let $\delta= \epsilon/\epsilon_r$. Then $\epsilon,\delta$ satisfy the hypotheses of \ref{smallradius}.

Let $(G,\mu)$ be an $\epsilon$-coherent massed graph; we will prove there is a versatile copy of $T$ in $G$. 
If $(G,\mu)$ is $(\delta,r+1)$-focussed, the result follows from \ref{smallradius}, so we assume not. Hence
there exists $Z\subseteq V(G)$ with $\mu(Z)\ge \delta$, such that $\mu(N^{r+1}_{G[Z]}[v])< \mu(Z)/2$ for each $v\in Z$.
Let $\mu'(X)=\mu(X)/\mu(Z)$ for all $X\subseteq Z$; then $(G[Z],\mu')$
is a massed graph. If it is $(\epsilon_r,r)$-coherent then the result follows from \ref{bigradius}, so we assume not,
for a contradiction. If there exist anticomplete subsets $A,B$ of $Z$ with $\mu'(A),\mu'(B)>\epsilon_r$, 
then $\mu(A), \mu(B)\ge \epsilon$, which is impossible. Thus there exists
$v\in Z$ such that $\mu'(N^r_{G[Z]}[v])\ge \epsilon_r$, and hence such that $\mu(N^r_{G[Z]}[v])\ge \epsilon_r \mu(Z)\ge \epsilon$.
But $\mu(N^{r+1}_{G[Z]}[v])< \mu(Z)/2$ from the choice of $Z$, and so $\mu(Z\setminus N^{r+1}_{G[Z]}[v])>\mu(Z)/2\ge \epsilon$,
a contradiction since the two sets $N^r_{G[Z]}[v]$ and $Z\setminus N^{r+1}_{G[Z]}[v]$ are anticomplete.
This proves \ref{catversion2}.~\bbox

\section{Parallels with $\chi$-boundedness}

An ideal is {\em $\chi$-bounded} if there is a function $f$ such that $\chi(G)\le f(\omega(G))$ for each graph $G$ in 
the ideal. Such ideals have been studied intensively, and it turns out that incoherent ideals and $\chi$-bounded ideals are in some ways very similar.
Here are two instances:
\begin{itemize}
\item 
Take a graph $H$, and let $\mathcal{I}$ be the ideal of all $H$-free graphs. The Gy\'arf\'as-Sumner 
conjecture~\cite{gyarfastree, sumner} asserts that $\mathcal{I}$ is $\chi$-bounded if and only if $H$ is a forest, and a conjecture
of~\cite{caterpillars} asserts that $\mathcal{I}$ is incoherent if and only if $H$ is a forest. The second is now 
(very recently) a theorem~\cite{treesanti}, and the ``only if'' part of the first is true, and the ``if'' part has been shown 
for some forests. 
\item A {\em hole} in $G$ means an induced cycle of length at least four. Let
$\mathcal{I}$ be the class of all graphs with no hole of length at least $k$, for some fixed integer $k$.
A theorem of~\cite{longholes} says that $\mathcal{I}$ is $\chi$-bounded, and a theorem of~\cite{bonamy} says that
$\mathcal{I}$ is incoherent.
\end{itemize}
But the parallel does not always work, and in fact neither property implies the other. 
Here are two examples showing this, one for each direction:
\begin{itemize}
\item The ideal of all perfect graphs is $\chi$-bounded, and indeed so is the ideal of all graphs
with no odd hole~\cite{oddholes}, but an example
of Fox~\cite{fox} shows that these ideals are coherent. 
\item Fix a graph $H$, and let $\mathcal{I}$
be the ideal of all graphs such that no induced subgraph is isomorphic to a subdivision of $H$. 
Then $\mathcal{I}$ is not necessarily $\chi$-bounded~\cite{sevenpoles}, but our main result
proves that it is incoherent.
\end{itemize}


There are a number of hard results and longstanding open conjectures about ideals that are not $\chi$-bounded, and it is entertaining 
to try their parallels
for coherent ideals. (The proofs below are just sketched.)

It is proved in~\cite{holeparity} that every ideal that is not $\chi$-bounded contains cycles of all 
lengths modulo $k$, for every integer $k\ge 1$. The same is not true for coherent ideals, as the example of 
Fox~\cite{fox} shows; a coherent ideal need not contain a cycle of odd length more than three, and in particular need not contain
a cycle of length 1 modulo 6. But it follows from \ref{mainthm} that for all integers $\ell$, every coherent ideal 
contains a cycle of length $2\ell$ modulo $k$, and hence contains one of 
every length modulo
$k$ if $k$ is odd. To see this, choose $\epsilon>0$
very small, and choose an $\epsilon$-coherent graph $G$ from the ideal. By \ref{mainthm}, $G$ contains a $P$-filleting of 
complete
graph $H$ of some large (constant) size, where $P$ is a Hamilton path of $H$. Choose many disjoint subpaths of $P$,
each of length $2k$, with no edges joining them. Let these paths be $P_1\ll P_n$ say, and let the $i$th vertex of $P_j$
be $v^i_j$. For each $i$, and $1\le h<k\le n$, there is a path of the $P$-filleting that joins $v^i_h, v^i_j$, say $Q_i^{h,j}$; 
and by 
Ramsey's theorem, we may choose many of the paths $P_j$ such that all the paths $Q_i^{h,j}$ have the same length modulo $k$
depending on $i$. (Redefine $n$, and renumber $P_1\ll P_n$ so that this holds.) Let $R_i$ be the union of 
$$Q_i^{1,2}, Q_i^{2,3}\ll Q_i^{k,k+1};$$
then $R_i$ has length divisible by $k$. But then for any $t$ modulo $k$, the union of $R_1, R_{t+1}$ and subpaths 
of $P_1$ and $P_{k+1}$ (both of length 
$t$) makes a hole of length $2t$ modulo $k$, and this cycle belongs to the ideal.

It is conjectured in~\cite{holeseq} that in every graph with huge chromatic number and bounded clique number, there
are $k$ holes with consecutive lengths. The same is not true in $\epsilon$-coherent graphs $G$ with $\epsilon$ very small,
because they need not have odd holes; but perhaps there must always be $k$ even holes with successive lengths differing by two?

It is proved in \cite{rainbow} that, in any colouring of a graph with huge chromatic number and bounded clique number, 
some induced $k$-vertex path is {\em rainbow} (that is, all its vertices have different colours), and no other types of 
connected subgraph have this property. What if we colour a graph which is $\epsilon$-coherent for $\epsilon$ very small?
Then we can do better than just paths; the results of this paper show that we can get a rainbow copy of any caterpillar.  
Each colour class has cardinality at most $2\epsilon|G|$, so by grouping the 
colour classes, we can partition the vertex set into many disjoint sets each of about the same size 
(differing by at most $2\epsilon |G|$), and each a union of colour classes. Then \ref{findcat} gives a 
copy $T'$ of $T$ with at most one vertex from each
block of the partition; and in particular, $T'$ is rainbow. Actually, we can do even better; results of 
~\cite{treesanti} show we can get a rainbow copy of any forest.

If we direct the edges of a graph with huge chromatic number and bounded clique number, some digraphs must be present
as induced subdigraphs. For instance, it is proved in \cite{orientations} that every oriented star has this property,
and so does a three-edge path where both ends point outwards. What if we direct the edges of a graph 
which is $\epsilon$-coherent for $\epsilon$ very small? Now much less is true. We need not get a directed two-edge path, 
because of
Fox's example from~\cite{fox} (this is a comparability graph, and so can be directed so that there is no induced 
directed two-edge path). We also need not get an outdirected 3-star (a tree with four vertices, three of them adjacent 
from the fourth). To see this, fix $\epsilon$, choose $k$ with $k\epsilon\ge 2$, and take $k$
disjoint sets $A_1\ll A_k$ each of the same size $n/k$ say, with $n$ large. 
Now take a random graph on $A_1\cup\cdots\cup A_k$ with average degree
$\log n$; with high probability the outcome is $\epsilon$-coherent, and its maximum degree is $O(\log(n))$. For
$1\le i\le k$ in turn, and for every pair of vertices $u,v\in A_{i+1}\cup\cdots\cup A_k$ with a common neighbour in $A_i$,
add an edge $uv$. Let the result be $G'$. Since this process is repeated only $k$ times and the maximum degree 
at most squares at each step,
the maximum degree of $G'$ is still less than $\epsilon n$. Now add more edges so that each $A_i$ is a clique, forming $G''$.
Thus $G''$ is $\epsilon$-coherent.
Orient every edge $uv$ of $G'$, from $u$ to $v$ if $u\in A_i$ and $v\in A_j$ where $i<j$, and arbitrarily if $u,v$ belong to the
same $A_i$. In the resulting digraph, there is no induced outdirected 3-star. (Incidentally, because of our 
main theorem \ref{mainthm} we always get a subdivision of $K_{2,3}$ as an induced subgraph, and however this 
is oriented it contains an outdirected 2-star; so we always get an outdirected 2-star.)



\begin{thebibliography}{99}
\bibitem{APPRS}
N. Alon, J. Pach, R. Pinchasi, R. Radoi\v ci\'c and M. Sharir, ''Crossing patterns of semi-algebraic
sets'', {\em J. Combinatorial Theory, Ser. A}, {\bf 111} (2005), 310--326.
\bibitem{bonamy} M. Bonamy, N. Bousquet, and S. Thomass\'e, ``The Erd\H{o}s-Hajnal conjecture for long holes and antiholes'',
{\em SIAM J. Discrete Math.}, {\bf 30} (2016), 1159--1164.
\bibitem{lagoutte} N. Bousquet, A. Lagoutte, and S. Thomass\'e, ``The Erd\H{o}s-Hajnal conjecture for paths and antipaths'', {\em J. Combinatorial Theory, Ser. B}, {\bf 113} (2015), 261--264.
\bibitem{hooks} K. Choromanski, D. Falik, A. Liebenau, V. Patel and M. Pilipczuk, ``Excluding hooks and their complements'', {\tt arXiv:1508.00634}.
\bibitem{maria} M. Chudnovsky, ``The Erd\H{o}s-Hajnal conjecture -- a survey'', {\em J. Graph Theory}
{\bf 75} (2014), 178--190.
\bibitem{longholes} M. Chudnovsky, A. Scott and P. Seymour, ``Induced subgraphs of graphs with large chromatic number.
III. Long holes'', {\em Combinatorica} {\bf 37} (2017), 1057--72.
\bibitem{orientations} M. Chudnovsky, A. Scott and P. Seymour, 
``Induced subgraphs of graphs with large chromatic number. XI.
Orientations'', {\em European Journal of Combinatorics} 76 (2019),
53--61,
{\tt arXiv:1711.07679}.
\bibitem{treesanti} M. Chudnovsky, A. Scott, P. Seymour and S. Spirkl,
``Pure pairs. I. Trees and linear anticomplete pairs'' 
(manuscript March 2018), {\tt arXiv:1809.00919}.
\bibitem{cs0} M. Chudnovsky and P. Seymour, 
``Excluding paths and antipaths'', {\em Combinatorica} {\bf 35} (2015), 389--412.
\bibitem{cz} M. Chudnovsky and Y. Zwols, ``Large cliques or stable sets in graphs with no four-edge path and no
five-edge path in the complement'', {\em J. Graph Theory}, {\bf 70} (2012), 449--472.
\bibitem{E0} P. Erd\H{o}s, ``Some remarks on the theory of graphs'', {\em Bull. Amer. Math. Soc.} {\bf 53} (1947), 292--294.
\bibitem{EH0} P. Erd\H{o}s and A. Hajnal, ``On spanned subgraphs of graphs'',
{\em Graphentheorie und Ihre Anwendungen} (Oberhof, 1977), \verb++{www.renyi.hu/\raisebox{-1ex}{\textasciitilde}p\_erdos/1977-19.pdf}.
\bibitem{EH}  P. Erd\H{o}s and A. Hajnal, ``Ramsey-type theorems'',
{\em  Discrete Applied Mathematics} {\bf 25} (1989), 37--52.
\bibitem{ES} P. Erd\H{o}s and G. Szekeres, ``A combinatorial problem in geometry'', {\em Compositio Mathematica}
{\bf 2} (1935), 463--470.
\bibitem{fox} J. Fox, ``A bipartite analogue of Dilworth's theorem'',
{\em Order} {\bf 23} (2006), 197--209.
\bibitem{fp}
J. Fox and J. Pach, ``Erd\H{o}s-Hajnal-type results on intersection patterns of geometric objects'', in
{\em Horizon of Combinatorics} (G.O.H. Katona et al., eds.), Bolyai Society Studies in Mathematics,
Springer, 79--103, 2008.
\bibitem{gyarfastree}
A. Gy\'arf\'as, ``On Ramsey covering-numbers'',
{\em Coll. Math. Soc. J\'anos Bolyai}, in {\em Infinite and Finite Sets},
North Holland/American Elsevier, New York (1975), 10.
\bibitem{caterpillars} A. Liebenau, M. Pilipczuk, P. Seymour and S. Spirkl, ``Caterpillars in Erd\H{o}s-Hajnal'',
{\em J. Combinatorial Theory, Ser. B}, 136 (2019), 33--43, {\tt arXiv:1810.00811}.
\bibitem{sevenpoles} A. Pawlik, J. Kozik, T. Krawczyk, M. Laso\'{n}, P. Micek, W. T. Trotter and B. Walczak,
``Triangle-free intersection graphs of line segments with large chromatic number'',
{\em J. Combinatorial Theory, Ser. B}, {\bf 105} (2014), 6--10.
\bibitem{rodl} V. R\"odl, ``On universality of graphs with uniformly distributed edges'',
{\em Discrete Math.} {\bf 59} (1986), 125--134.
\bibitem{oddholes} A. Scott and P. Seymour, 
``Induced subgraphs of graphs with large chromatic number.  I. Odd holes'', 
{\em J. Combinatorial Theory, Ser. B}, {\bf 121} (2016), 68--84.
\bibitem{holeseq} A. Scott and P. Seymour,  ``Induced subgraphs of graphs with large chromatic
number. IV. Consecutive holes'', {\em J. Combinatorial Theory, Ser. B},
132 (2018), 180--235, {\tt arXiv:1509.06563}.
\bibitem{rainbow} A. Scott and P. Seymour, 
``Induced subgraphs of graphs with large chromatic number. IX.
Rainbow paths'', {\em Electronic J. Combinatorics}, 24.2:\#P2.53, 2017.
\bibitem{holeparity} A. Scott and P. Seymour,
``Induced subgraphs of graphs with large chromatic number. X. Holes with
specific residue'', submitted for publication, {\tt arXiv:1705.04609}.
\bibitem{sumner}
D.P. Sumner, ``Subtrees of a graph and chromatic number'', in
{\em The Theory and Applications of Graphs}, (G. Chartrand, ed.),
John Wiley \& Sons, New York (1981), 557--576.
\end{thebibliography}
\end{document}